\documentclass{article}
\usepackage{amsmath}
\usepackage{amssymb}
\usepackage{graphicx}
\usepackage{float}
\usepackage{cite}
\usepackage{subfigure}
\usepackage{pgfplotstable}
\usepackage{booktabs}
\usepackage{caption}
\usepackage{subfig}
\usepackage[top=1.5in, bottom=0.5in, left=1.5in, right=1in]{geometry}
\begin{document}
\begin{center}
\textbf{\Large{Determination of a space-dependent force function in the one-dimensional wave equation}}
\end{center}
S.O. Hussein and D. Lesnic 
\\
\textit{Department of Applied Mathematics, University of Leeds, Leeds LS2 9JT, UK}\\
E-mails: ml10soh@leeds.ac.uk, D.Lesnic@leeds.ac.uk\\
\\
\textbf{\large Abstract.}
The determination of an unknown spacewice dependent force function acting on a vibrating string from over-specified Cauchy boundary data is investigated numerically using the boundary element method (BEM) combined with a regularized method of separating variables.
This linear inverse problem is ill-posed since small errors in the input data cause large errors in the output force solution. Consequently, when the input data is contaminated with noise we use the Tikhonov regularization method in order to obtain a stable solution. The choice of the regularization parameter is based on the L-curve method. Numerical results show that the solution is accurate for exact data and stable for noisy data.

\section{Introduction}
The wave equation governs many physical problems such as the vibrations of a spring or membrane, acoustic scattering, etc. When it comes to mathematical modeling probably the most investigated are the direct and inverse acoustic scattering problems, see e.g. \cite{kresscoltan2013}.

On the other hand, inverse source/force problems for the wave equation have been less investigated. It is the objective of this study to investigate such an inverse force problem for the hyperbolic wave equation. The initial attempt is performed for the case of a one-dimensional vibrating string, but we have in mind extensions to higher dimensions in an immediate future work. The forcing function is assumed to depend only upon the single space variable in order to ensure uniqueness of the solution. The theoretical basis for our numerical investigation is given in \cite{candunn70} where the uniqueness of solution of the inverse spacewise dependent force function for the one-dimensional wave equation has been established. The authors of \cite{candunn70} have also given conditions to be satisfied by the force function in order to ensure continuous dependence upon the data and furthermore, they proposed two methods based on linear programming and the least-squares method. However, no numerical results were presented and it is the main purpose of our study to develop an efficient numerical solution for this inverse linear, but ill-posed problem.

Because the wave speed is assumed constant, the most suitable numerical method for discretising the wave equation in this case is the boundary element method (BEM), see [1-3]. Moreover, because an inhomogeneous source/force term is present in the governing equation, it is convenient to exploit the linearity of the problem by applying the principle of superposition. This recasts into splitting the original problem into a direct problem with no force, and an inverse problem with force, but with homogeneous boundary and initial conditions. This is explained in section 2 where the mathematical formulation of the inverse problem under investigation is also given. Whilst the former problem requires a numerical solution such as the BEM, as described in Section 3, the latter problem is ameanable to a separation of variables series solution with unknown coefficients. Upon truncating this series, the problem recasts as an ordinary linear least-squares problem which has to be regularized since the resulting system of linear equations is ill-conditioned, the original problem being ill-posed. The choice of the regularization parameter introduced by this technique is important for the stability of the numerical solution and in our study this is based on the heuristic L-curve criterion,\cite{hansen2001}. All this latter analysis is described in detail in Section 4. Numerical results are illustrated and discussed in Sections 5 and 6 and conclusions and future work are provided in Section 7. 
\section{Mathematical Formulation} 
The governing equation for a vibrating string of length $L>0$ acted upon a space-dependent force $f(x)$ is given by the one-dimensional wave equation 
\begin{equation}
u_{tt}=c^2 u_{xx}+f(x), \ \ \ \ \ \ \ \ \ \ \   x\in(0,L)\times (0,\infty), \label{eq1}         
\end{equation}
where $u$ represents the displacement and $c>0$ is the speed of sound.

Equation \eqref{eq1} has to be solved subject to the initial conditions
\begin{equation}
u(x,0)=u_0(x), \ \ \ \ \ \ \ \ x\in[0,L],   \label{eq2}
\end{equation}
\begin{equation}
 u_{t}(x,0)=v_{0}(x), \ \ \ \ \ \ \  x\in[0,L], \label{eq3}
\end{equation}
where $u_0$ and $v_{0}$ represent the initial displacement and velocity, respectively, and to the Dirichlet boundary conditions 
\begin{eqnarray}
u(0,t)=p_{0}(t), \ \ \ \ \ \ \ \ \ \ \ \ \ \ \ \ \ \ \ \ \ \ \ \ \ \ \ \ \ \ t\in[0,\infty), \label{eq4}
\end{eqnarray}
\begin{eqnarray}
\mu u(L,t)+(1-\mu)u_{x}(L,t)=p_{L}(t), \ \ \ t\in [0,\infty), \label{eq5}
\end{eqnarray}
where $\mu \in \lbrace 0,1 \rbrace$ with $\mu =1$ for the Dirichlet boundary condition and $\mu =0$ for the Neumann boundary condition. In \eqref{eq4} and \eqref{eq5}, $p_{0}$ and $p_{L}$ are given functions satisfying the compatibility conditions
\begin{eqnarray}
p_{0}(0)=u_{0}(0), \ \ \ p_{L}(0)=\mu u_{0}(L)+(1-\mu)u^{\prime}_{0}(L). \label{eq6}
\end{eqnarray} 

If the force $f$ is given, then equations \eqref{eq1}-\eqref{eq6}
form a direct well-posed problem for $u(x,t)$ which can be solved using the BEM for example, \cite{benmansour93}. However, if the force function $f$ is unknown then clearly the above equations are not sufficient to determine the pair solution $(u(x,t),f(x))$. Then, as suggested in \cite{candunn70}, we supply the above system of equations with the measurement of the flux tension of the string at the end $x=0$, namely  
\begin{equation}
u_{x}(0,t)=q_{0}(t), \ \ \ t\in[0,T], \label{eq117}
\end{equation}
where $q_{0}$ is a given function over a time of interest $T>0$.
Then the inverse problem under investigation requires determining the pair solution$(u(x,t),f(x))$ satisfying equations \eqref{eq1}-\eqref{eq117}. Remark that we have to restrict $f$ to depend on $x$ only since otherwise, if $f$ depends on both $x$ and $t$, we can always add to $u(x,t)$ any function of the form $t^{2}x^{2}(x-L)^{2}U(x,t)$ with arbitrary  $U\in C^{2,1}([0,L]\times[0,\infty))$  and still obtain another solution satisfying \eqref{eq1}-\eqref{eq117}. Note that the unknown force $f(x)$ depends on the space variable $x$, whilst the additional measurement \eqref{eq117} of the flux $q_{0}(t)$ depends on the time variable $t$.

It has been shown in \cite{candunn70} that the problem \eqref{eq1}-\eqref{eq117} has at most one solution, i.e. the uniqueness holds. Moreover, the solution depends continuously on the input data if $f \in C^{2}[0,L]$ with $supp(f)\subset(0,L)$ and
\begin{equation}
 max\lbrace |f(x)|, |f^{\prime}(x)|, |f^{\prime \prime}(x)| ; \ \  x\in [0,L] \rbrace \leq K,
\label{eq118}
\end{equation}
where $K$ is a known positive constant.

Due to the linearity of the inverse problem \eqref{eq1}-\eqref{eq117} it is convenient to split it into the form, \cite{candunn70},
\begin{equation}
u=v+w, \label{eq119}
\end{equation}
where $v$ satisfies the well-posed direct problem
\begin{eqnarray}
v_{tt}&=&c^{2}v_{xx}, \ \ \ (x,t)\in(0,L)\times(0,\infty),\label{eq1110}\\
v(x,0)&=&u_{0}(x), \ \ \  x\in[0,L],    \label{eq11v}\\
v_{t}(x,0)&=&v_{0}(x), \ \ \ x\in [0,L], \label{eq12v}\\
v(0,t)&=&p_{0}(t), \ \ \ t\in[0,\infty),   \label{eq13v}\\
\mu v(L,t)+(1-\mu) v_{x}(L,t)&=&p_{L}(t),  \ \ \ t\in[0,\infty). 
\label{eq14v}
\end{eqnarray}
and $(w,f)$ satisfies the ill-posed inverse problem
\begin{eqnarray}
w_{tt}&=&c^{2}w_{xx}+f(x), \ \ \ (x,t)\in (0,L)\times(0,\infty), \label{eq15w}\\
w(x,0)&=&w_{t}(x,0)=0, \ \ \ x\in[0,L], \label{eq16w}\\
w(0,t)&=&0, \ \ \ t\in[0,\infty), \label{eq17w}\\
\mu w(L,t)+(1-\mu) w_{x}(L,t)&=&0, \ \ \ t\in[0,\infty), \label{eq18w}\\
w_{x}(0,t)&=&q_{0}(t)-v_{x}(0,t), \ \ \ t\in[0,T]. \label{eq19w}
\end{eqnarray}
The problems \eqref{eq1}-\eqref{eq117}, \eqref{eq1110}-\eqref{eq14v} and \eqref{eq15w}-\eqref{eq19w} are schematically depicted in Figure 1.

Observe that we could also control the Dirichlet data \eqref{eq4} instead of the Neumann data \eqref{eq117}. We remark that the solution of the direct and well-posed problem \eqref{eq1110}-\eqref{eq14v} has to be found numerically, say using the BEM, as described in the next section.

\begin{figure}[H]
\centering
\includegraphics[width=12.5cm]{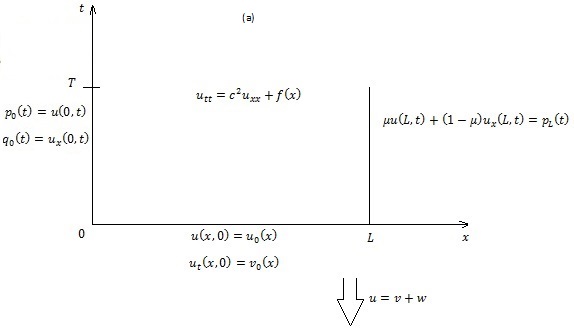}
\includegraphics[width=12.5cm]{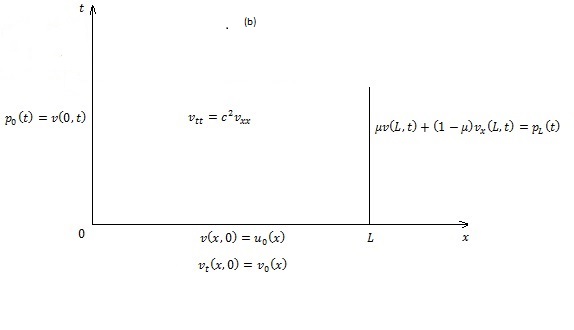}
\end{figure}
\newpage
\begin{figure}[H]
\centering
\includegraphics[width=12.5cm]{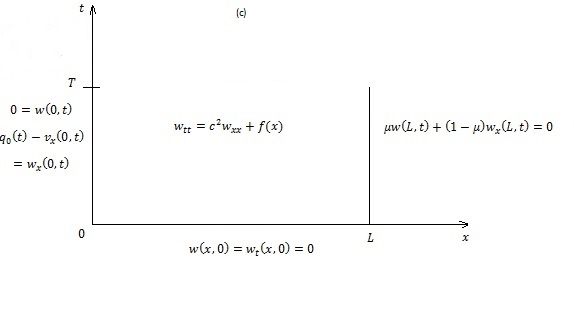}
\caption{Schematic diagram depicting the problems (a) \eqref{eq1}-\eqref{eq117}, (b) \eqref{eq1110}-\eqref{eq14v}, and (c) \eqref{eq15w}-\eqref{eq19w}.}
\label{fig:fg3}
\end{figure}

\section{The Boundary Element Method (BEM) for Solving the Direct Problem \textbf{\eqref{eq1110}}-\textbf{\eqref{eq14v}}}
The BEM for the one-dimensional wave equation \eqref{eq1110} is based on the application of integration by parts and the use of the fundamental solution, ([7], p.893),
\begin{equation}
u^{*}(x,t;\xi,\tau)=-\frac{1}{2c}H(c(t-\tau)-|x-\xi|),  \label{eq3.1}
\end{equation}
where $H$ is the Heaviside function. This results in the following boundary integral equation, \cite{benmansour93},
\begin{eqnarray}
2v(\xi,t)=v(\xi-ct,0)+v(\xi+ct,0)+\frac{1}{c}\int_{\xi-ct}^{\xi+ct}v_{t}(x,0) dx+v(L,t-(L-\xi)/c)\notag \\
+c\int_{0}^{t-(L-\xi)/c} v_{x}(L,\tau)d\tau+v(0,t-\xi/c)-c\int_{0}^{t-\xi/c} v_{x}(0,\tau)d\tau, \notag \\ 
(\xi,t)\in(0,L)\times (0,\infty). \label{eq3.2}
\end{eqnarray}
Equation \eqref{eq3.2} is valid if 
\begin{equation}
v(0,0)=v(L,0)=0, \ \ \ i.e.\ \ \ u_{0}(0)=u_{0}(L)=0. \label{eq3.3}
\end{equation} 
Otherwise, if this condition is not satisfied then we can work with the modified function
\begin{equation}
\widetilde{v}(x,t)=v(x,t)-\frac{u_{0}(L)-u_{0}(0)}{L}x-u_{0}(0) \label{eq3.4}
\end{equation}
which satisfies $\widetilde{v}(0,0)=\widetilde{v}(L,0)=0$.

It is very important to remark that in expression \eqref{eq3.2} the time and space coordinates must be within the domain $[0,L]\times [0,\infty)$ and the integrals must have their lower limit of integration smaller than the upper one. If any of these conditions are not satisfied the integrals are taken to be zero.

Equation \eqref{eq3.2} yields the interior solution $v(\xi,t)$ for $(\xi,t)\in (0,L)\times (0,\infty)$ of the wave equation \eqref{eq1110} in terms of the initial and boundary data. In general, at a boundary point only one Dirichlet, Neumann or Robin boundary condition is imposed and the first step of the BEM methodology requires the evaluation of the missing (unspecified) boundary data. For this, we need first to evaluate the boundary integral equation \eqref{eq3.2} at the end points $\xi \in \lbrace 0,L \rbrace$. A careful limiting process yields, \cite{benmansour93},
\begin{eqnarray}
v(0,t)&=&v(ct,0)+\frac{1}{c}\int_{0}^{ct}v_{t}(x,0)dx+v(L,t-L/c) \notag \\
& &+c{\bigg[}\int_{0}^{t-L/c}v_{x}(L,\tau)d\tau-\int_{0}^{t}v_{x}(0,\tau)d\tau{\bigg]}, \quad t\in (0,\infty),  \label{eq3.5a} \\
v(L,t)&=&v(0,t-L/c)+v(L-ct,0)+\frac{1}{c}\int_{L-ct}^{L}v_{t}(x,0)dx \notag \\
& &+c{\bigg[}\int_{0}^{t}v_{x}(L,\tau)d\tau-\int_{0}^{t-L/c} v_{x}(0,\tau)d\tau{\bigg]}, \quad t\in (0,\infty). \label{eq3.5b}
\end{eqnarray}
These equations also hold under the assumption \eqref{eq3.3}.

Since we want to calculate $v_{x}(0,t)$ only for $t\in[0,T]$, let us restrict the boundary integral equations \eqref{eq3.5a} and \eqref{eq3.5b} to the time interval $[0,T]$. For the numerical discretisation of the boundary integral equations \eqref{eq3.5a} and \eqref{eq3.5b} we divide the time interval $[0,T]$ into a series of $N$ small boundary elements $[t_{j-1},t_{j}]$ for $j=\overline{1,N}$, where for a uniform discretisation $t_{j}=j T/N$ for $j=\overline{0,N}$. Similarly, we divide the space interval $[0,L]$ into a series of $M$ small cells $[x_{i-1},x_{i}]$ for $i=\overline{1,M}$, where for a uniform discretisation $x_{i}=i L/M$ for $i=\overline{0,M}$. We then approximate the boundary and initial values as 
\begin{eqnarray}
v(0,\tau)=\sum_{j=1}^{N}\phi^{j}(\tau) v{^{0}_{j}},\ \  
v(L,\tau)=\sum_{j=1}^{N}\phi^{j}(\tau) v{^{L}_{j}}, \ \ \  \tau \in [0,T], \label{eq3.6a} \\
v_{x}(0,\tau)=\sum_{j=1}^{N}\theta^{j}(\tau) v{^{\prime 0}_{j}}, \ \ 
v_{x}(L,\tau)=\sum_{j=1}^{N}\theta^{j}(\tau) v{^{\prime L}_{j}}, \ \ \  \tau \in [0,T], \label{eq3.6b}\\
v(x,0)=\sum_{i=1}^{M}\psi_{i}(x) v{^{i}_{0}}, \ \ v_{t}(x,0)=\sum_{i=1}^{M}\psi_{i}(x) v{^{i}_{0}}, \ \ \ x\in [0,L], \label{eq3.6c}
\end{eqnarray}
where 
\begin{eqnarray}
v{^{0}_{j}}:=v(0,t_{j}), \ \ v{^{L}_{j}}:=v(L,t_{j}), \ \ v{^{\prime 0}_{j}}:=v_{x}(0,t_{j}), \ \ v{^{\prime L}_{j}}:=v_{x}(L,t_{j}), \ \ j=\overline{1,N}, \label{eq3.7a}\\
u{^{i}_{0}}:=v(x_{i},0), \ \ \ v{^{i}_{0}}:=v_{t}(x_{i},0), \ \ \ i=\overline{1,M}. \quad \quad \quad \quad \quad \label{eq3.7b}
\end{eqnarray} 
The functions $\phi^{j}$, $\theta^{j}$ and $\psi_{i}$ are interpolant, e.g. piecewise polynomial, functions chosen such that $\phi^{j}(t_{n})=\theta^{j}(t_{n})=\delta_{jn}$ for $j,n=\overline{1,N}$, $\psi_{i}(x_{m})=\delta_{im}$ for $i,m=\overline{1,M}$ , where $\delta_{jn}$ is the Kronecker delta symbol. For example, if $\theta^{j}(\tau)$ is a piecewise constant function then
\begin{equation}
\theta^{j} (\tau)=\chi_{(t_{j-1},t_{j}]}(\tau)=
\begin{cases}
1 \ \ \ \ \ \ \text{if\ \ \ $t \in(t_{j-1},t_{j}]$},
\\
0 \ \ \ \ \ \ \text{otherwise}, \label{eq3.22}
\end{cases}
\end{equation} 
where $\chi_{(t_{j-1},t_{j}]}$ represents the characteristic function of the interval $(t_{j-1},t_{j}]$. Thus $v_{x}(0,\tau)=v{^{\prime 0}_{j}}$ for $\tau \in(t_{j-1},t_{j}]$, etc. We also have that $\int_{0}^{t_{n}}\theta^{j}(\tau) d\tau=t_{j}-t_{j-1}$\ \ \  for \ \ \ $j=\overline{1,n}$. 

Using the approximations \eqref{eq3.6a}-\eqref{eq3.7b} into the equations \eqref{eq3.5a} and \eqref{eq3.5b} we obtain, for $n=\overline{1,N}$,
\begin{eqnarray}
& &v_{n}^{0}+cv_{n}^{\prime 0}\int_{0}^{t_{n}}\theta ^{n}(\tau) d\tau-\phi^{n}(t_{n}-L/c)v_{n}^{L}-cv_{n}^{\prime L}\int_{0}^{t_{n}-L/c}\theta^{n}(\tau)d\tau \notag \\
& &=\sum_{j=1}^{n-1}\phi^{j}(t_{n}-L/c)v_{j}^{L}+\sum_{i=1}^{M}\psi_{i}(ct_{n})u_{0}^{i}+c\sum_{j=1}^{n-1}v_{j}^{\prime L}\int_{0}^{t_{n}-L/c}\theta^{j}(\tau) d\tau \notag \\
& &- c\sum_{j=1}^{n-1}v_{j}^{\prime 0}\int_{0}^{t_{n}}\theta^{j}(\tau) d\tau+\frac{1}{c}\sum_{i=1}^{M}v{^{i}_{0}}\int_{0}^{ct_{n}}\psi_{i}(x) dx=:F \label{eq3.8a}
\end{eqnarray}
\\ and
\begin{eqnarray}
& &v_{n}^{L}-cv_{n}^{\prime L}\int_{0}^{t_{n}}\theta ^{n}(\tau) d\tau-\phi^{n}(t_{n}-L/c)v_{n}^{0}+cv_{n}^{\prime 0}\int_{0}^{t_{n}-L/c}\theta^{n}(\tau)d\tau \notag \\
& &=\sum_{j=1}^{n-1}\phi^{j}(t_{n}-L/c)v_{j}^{0}+\sum_{i=1}^{M}\psi_{i}(L-ct_{n})u_{0}^{i}+c\sum_{j=1}^{n-1}v_{j}^{\prime L}\int_{0}^{t_{n}}\theta^{j}(\tau) d\tau \notag \\
& &- c\sum_{j=1}^{n-1}v_{j}^{\prime 0}\int_{0}^{t_{n}-L/c}\theta^{j}(\tau) d\tau+\frac{1}{c}\sum_{i=1}^{M}v{^{i}_{0}}\int_{L-ct_{n}}^{L}\psi_{i}(x) dx=:G. \label{eq3.8b}
\end{eqnarray}
Denoting
\begin{eqnarray}
A=c\int_{0}^{t_{n}}\theta^{n}(\tau)d\tau,\ \ \ 
B=\phi^{n}(t_{n}-L/c),\ \ \ 
D=c \int_{0}^{t_{n}-L/c} \theta^{n}(\tau) d\tau,  \label{eq3.9}
\end{eqnarray}
equations \eqref{eq3.8a} and \eqref{eq3.8b} can be rewritten as
\begin{eqnarray}
& &v_{n}^{0}+Av_{n}^{\prime 0}-Bv_{n}^{L}-Dv_{n}^{\prime L}=F, \label{eq3.10a}\\
& &v_{n}^{L}-Av_{n}^{\prime L}-Bv_{n}^{0}+Dv_{n}^{\prime 0}=G. \label{eq3.10b}
\end{eqnarray}
At each time $t_{n}$ for $n=\overline{1,N}$, the system of equations \eqref{eq3.10a} and \eqref{eq3.10b} represents a time marching BEM technique in which the values of $F$ and $G$ are expressed in terms of the previous values of the solution at the times $t_{1},...,t_{n-1}$. Note that upon the imposition of the initial conditions \eqref{eq11v} and \eqref{eq12v} we know 
\begin{eqnarray}
u_{0}^{i}=v(x_{i},0)=u_{0}(x_{i}), \ \ \ v_{0}^{i}=v_{t}(x_{i},0)=v_{0}(x_{i}), \ \ \ i=\overline{1,M}. \label{eq3.11}
\end{eqnarray}

The system of equations \eqref{eq3.10a} and \eqref{eq3.10b} contains 2 equations with 4 unknowns. Two more equations are known from the boundary conditions \eqref{eq13v} and \eqref{eq14v}, namely
\begin{eqnarray}
v_{n}^{0}=v(0,t_{n})=p_{0}(t_{n})=:p_{0}^{n}, \ \ \  n=\overline{1,N}, \ \ \ \ \ \ \ \ \ \ \ \ \ \ \ \ \ \ \ \ \  \ \ \ \ \ \ \ \ \ \ \ \  \ \ \ \ \ \ \ \ \ \ \ \ \label{eq3.12a}\\
\mu v_{n}^{L}+(1-\mu)v_{n}^{\prime L}=\mu v(L,t_{n})+(1-\mu)v_{x}(L,t_{n})=p_{L}(t_{n})=:p_{L}^{n}, \ \ \ n=\overline{1,N}. \label{eq3.12b}
\end{eqnarray}

The solution of the system of equations \eqref{eq3.10a},   \eqref{eq3.10b}, \eqref{eq3.12a} and \eqref{eq3.12b} can be expressed explicitly at each time step $t_{n}$ for $n=\overline{1,N}$ and is given by: \\
(a) For $\mu =1$, i.e. the Dirichlet problem \eqref{eq1110}-\eqref{eq14v} in which equation \eqref{eq14v} is given by 
\begin{eqnarray}
v(L,t)=p_{L}(t), \ \ \  t\in [0,\infty),  \label{eq3.13}
\end{eqnarray}
and equation \eqref{eq3.12b} yields
\begin{eqnarray}
v_{n}^{L}=p_{L}^{n},  \ \ \ n=\overline{1,N}, \label{eq3.14}
\end{eqnarray}
the unspecified boundary values are the Neumann flux values. Introduction of \eqref{eq3.12a} and \eqref{eq3.14} into \eqref{eq3.10a} and \eqref{eq3.10b} yields the simplified system of two equations with two unknowns given by 
\begin{eqnarray}
& &Av_{n}^{\prime 0}-Dv_{n}^{\prime L}=F-p_{0}^{n}+Bp_{L}^{n}=:\widetilde{F}, \label{eq3.15a}\\
& &Dv_{n}^{\prime 0}-Av_{n}^{\prime L}=G-p_{L}^{n}+Bp_{0}^{n}=:\widetilde{G}. \label{eq3.15b}
\end{eqnarray}
Application of Cramer's rule immediately yields the solution
\begin{eqnarray}
v_{n}^{\prime 0}=\frac{D\widetilde{G}-A\widetilde{F}}{D^{2}-A^{2}}, \ \ \  v_{n}^{\prime L}=\frac{A\widetilde{G}-D\widetilde{F}}{D^{2}-A^{2}}. \label{eq3.16}
\end{eqnarray}
(b) For $\mu=0$, i.e. the mixed problem \eqref{eq1110}-\eqref{eq14v} in which equation \eqref{eq14v} is given by 
\begin{eqnarray}
v_{x}(L,t)=p_{L}, \ \ \ t\in[0,\infty) \label{eq3.17}
\end{eqnarray}
and equation \eqref{eq3.12b} yields
\begin{eqnarray}
v_{n}^{\prime L}=p_{L}^{n}, \ \ \ n=\overline{1,N}, \label{eq3.18}
\end{eqnarray}
the unspecified boundary values are the Neumann data at $x=0$ and the Dirichlet data at $x=L$. Introduction of \eqref{eq3.12a} and \eqref{eq3.18} into \eqref{eq3.10a} and \eqref{eq3.10b} yields
\begin{eqnarray}
& &Av_{n}^{\prime 0}-Bv_{n}^{L}=F-p_{0}^{n}+Dp_{L}^{n}=:\widetilde{\widetilde{F}}, \label{eq3.19a}\\
& &Dv_{n}^{\prime 0}+v_{n}^{L}=G+Bp_{0}^{n}+Ap_{L}^{n}=: \widetilde{\widetilde{G}}. \label{eq3.19b}
\end{eqnarray}
This yields the solution
\begin{eqnarray}
v_{n}^{\prime 0}=\frac{\widetilde{\widetilde{F}}+B\widetilde{\widetilde{G}}}{A+DB}, \ \ \  v_{n}^{L}=\frac{A\widetilde{\widetilde{G}}-D\widetilde{\widetilde{F}}}{A+DB}. \label{eq3.20}
\end{eqnarray}

Alternatively, instead of employing a time-marching BEM it is also possible to employ a global BEM by assembling \eqref{eq3.8a} and \eqref{eq3.8b} as a full system of $2N$ linear equations with $4N$ unknown $v_{j}^{0}, v_{j}^{L}, v_{j}^{\prime 0}, v_{j}^{\prime L} $ for $j=\overline{1,N}$, namely,
\begin{eqnarray}
& &\sum_{j=1}^{n}{\bigg[}cv{^{\prime 0}_{j}}\int_{0}^{t_{n}}\theta ^{j}(\tau) d\tau-cv{^{\prime L}_{j}}\int_{0}^{t_{n}-L/c}\theta^{j}(\tau)d\tau-\phi^{j}(t_{n}-L/c)v{^{L}_{j}}{\bigg]} \notag \\
& &=\sum_{i=1}^{M}{\bigg[}\psi_{i}(ct_{n})u{^{i}_{0}}+\frac{1}{c}v{^{i}_{0}}\int_{0}^{ct_{n}}\psi_{i}(x) dx{\bigg]}-v_{n}^{0}, \quad n=\overline{1,N}, \label{eq3.21a}
\end{eqnarray}
\\ and
\begin{eqnarray}
& &\sum_{j=1}^{n}{\bigg[}cv{^{\prime L}_{j}}\int_{0}^{t_{n}}\theta ^{j}(\tau) d\tau-cv{^{\prime 0}_{j}}\int_{0}^{t_{n}-L/c}\theta^{j}(\tau)d\tau+\phi^{j}(t_{n}-L/c)v{^{0}_{j}}{\bigg]} \notag \\
& &=-\sum_{i=1}^{M}{\bigg[}\psi_{i}(L-ct_{n})u{^{i}_{0}}+\frac{1}{c}v{^{i}_{0}}\int_{L-ct_{n}}^{L}\psi_{i}(x) dx{\bigg]}+v_{n}^{L} \quad n=\overline{1,N}. \label{eq3.21b}
\end{eqnarray}
The other $2N$ equations are given by \eqref{eq3.12a} and \eqref{eq3.12b}. Introduction of \eqref{eq3.12a} and \eqref{eq3.12b} into \eqref{eq3.21a} and \eqref{eq3.21b} finally results in a linear system of $2N$ algebraic equations with $2N$ unknowns which can be solved using a Gaussian elimination procedure.

Once all the boundary values been determined accurately, the interior solution can be obtained explicitly using equation \eqref{eq3.2}. This gives 
 \begin{eqnarray}
2v(\xi,t_{n})=\sum_{j=1}^{n}{\bigg[}\phi^{j}(t_{n}-(L-\xi)/c)v^{L}_{j}+\phi^{j}(t_{n}-\xi/c)v^{0}_{j}{\bigg]}\quad \quad \quad \quad \quad \quad \quad \notag \ \ \\
+c\sum_{j=1}^{n}{\bigg[}v^{\prime L}_{j}\int_{0}^{t_{n}-(L-\xi)/c}\theta^{j}(\tau)d\tau 
-v^{\prime 0}_{j}\int_{0}^{t_{n}-\xi/c}\theta^{j}(\tau)d\tau{\bigg]}
+\sum_{i=1}^{M}{\bigg[}\psi_{i}(\xi-ct_{n}) \notag \\
+\psi_{i}(\xi+ct_{n}){\bigg]}u^{i}_{0}
+\frac{1}{c}\sum_{i=1}^{M}v^{i}_{0}\int_{\xi-ct_{n}}^{\xi+ct_{n}}\psi_{i}(x)dx,  \quad \quad \quad n=\overline{1,N}, \quad \xi \in(0,1). \ \ \ \label{eq3.23} 
\end{eqnarray}
In \eqref{eq3.23}, for the piecewise constant interpolation \eqref{eq3.22}, 
\begin{eqnarray*}
\int_{0}^{t_{n}-\xi/c}\theta^{j}(\tau)d\tau&=& H(t_{n}-t_{j-1}-\xi/c)(t_{j}-t_{j-1}), \\
\int_{0}^{t_{n}-(L-\xi)/c}\theta^{j}(\tau)d\tau&=& H(t_{n}-t_{j-1}-(L-\xi)/c)(t_{j}-t_{j-1}).
\end{eqnarray*} 

As experienced in \cite{benmansour93}, by comparing the BEM with the method of characteristics one observes that by choosing  $\Delta t=T/N$ and $\Delta x=L/M$ such that the Courant number $c\Delta t/\Delta x=cTM/(NL)=1$ then, the grid is set to fit the characteristic net and the error involved is made arbitrarily small by reducing interpolation between points.

The flux $v_{x}(0,t)$ obtained numerically using the BEM is then introduced into \eqref{eq19w} and the inverse problem \eqref{eq15w}-\eqref{eq19w} for the pair solution $(w(x,t),f(x))$ is solved using the method described in the next section.
\section{Method for solving the inverse problem \eqref{eq15w}-\eqref{eq19w}}
The method employed for solving the inverse problem \eqref{eq15w}-\eqref{eq19w} is based on the separation of variables which gives that an approximate solution to the problem \eqref{eq15w}-\eqref{eq19w} is given by, \cite{candunn70},
\begin{eqnarray}
w_{K}(x,t;\underline{b})=\frac{\sqrt{2}}{c^{2}}\sum_{k=1}^{K}\frac{b_{k}}{\lambda^{2}_{k}}(1-\cos(c\lambda_{k} t))\sin(\lambda_{k} x), \ \ \ (x,t)\in[0,L]\times[0,\infty), \label{eq4.1}
\end{eqnarray}
\begin{eqnarray}
f_{K}(x)=\sqrt{2}\sum_{k=1}^{K} b_{k}\sin(\lambda_{k} x), \ \ \  x\in(0,L),  \label{eq4.2}
\end{eqnarray}
where $K$ is a truncation number and 
 \begin{equation}
\lambda_{k} =
\begin{cases}
\frac{k\pi}{L} \ \ \ \ \ \ \ \ \ \ \text{if} \ \ \ \ \ \ \ \ \ \ \ \mu=1,
\\
\\ \frac{(k-\frac{1}{2})\pi}{L} \ \ \ \ \ \ \text{if} \ \ \ \ \ \ \ \ \ \ \ \mu=0.    \label{eq4.3}
\end{cases}
\end{equation}

The coefficients $\underline{b}=(b_{k})_{k=\overline{1,K}}$ are to be determined by imposing the additional boundary condition \eqref{eq19w}. This results in
\begin{eqnarray}
q_{0}(t)-v_{x}(0,t)=:g(t)=\frac{\partial{w_{k}}}{\partial{x}}(0,t;\underline{b})=\frac{\sqrt{2}}{c^{2}}\sum_{k=1}^{K}\frac{b_{k}}{\lambda_{k}}(1-\cos(c\lambda_{k} t)), \ \ \ t\in[0,T]. \label{eq4.4}
\end{eqnarray}

In practice, the additional observation \eqref{eq117} comes from measurement which is inherently contaminated with errors. We therefore model this by replacing the exact data $q_{0}(t)$ by the noisy data 
\begin{eqnarray}
q_{0}^{\epsilon}(t_{n})=q_{0}(t_{n})+\epsilon, \ \ \ n=\overline{1,N},\label{eq4.5}
\end{eqnarray}
where $\epsilon$ are $N$ random noisy variables generated (using the Fortran NAG routine G05DDF) from a Gaussian normal distribution with mean zero and standard deviation $\sigma$ given by 
\begin{eqnarray}
\sigma=p\% \times max_{t\in[0,T]}\left|q_{0}(t)\right|, \label{eq4.6}
\end{eqnarray}
where $p\%$ represents the percentage of noise. The noisy data \eqref{eq4.5} also induces noise in $g$ as given by
\begin{eqnarray}
g^{\epsilon}(t_{n})=q^{\epsilon}_{0}(t_{n})-v_{x}(0,t_{n})=g(t_{n})+\epsilon, \ \ \ n=\overline{1,N}. \label{eq4.7}
\end{eqnarray}

Then we apply the condition \eqref{eq4.4} with $g$ replaced by $g^{\epsilon}$ in a least-squares penalised sense by minimizing the Tikhonov functional
\begin{eqnarray}
\mathcal{J}(\underline{b}):=\sum_{n=1}^{N}{\bigg[}\frac{\sqrt{2}}{c^{2}}\sum_{k=1}^{K}\frac{b_{k}}{\lambda_{k}}(1-\cos(c\lambda_{k} t_{n}))-g^{\epsilon}(t_{n}){\bigg]}^{2}+\lambda \sum_{k=1}^{K}b_{k}^{2}, \label{eq4.8}
\end{eqnarray}
where $\lambda \geq 0$ is a regularization parameter to be prescribed according to some criterion, e.g. the L-curve criterion, \cite{hansen2001}.

Denoting
\begin{eqnarray}
\underline{g}^{\epsilon}=(g^{\epsilon}(t_{n}))_{n=\overline{1,N}}, \ \ \ Q_{nk}=\frac{\sqrt{2}(1-\cos(c\lambda_{k} t_{n}))}{c^{2}\lambda_{k}}, \ \ \ n=\overline{1,N}, \ \ \ k=\overline{1,K}, \label{eq4.9}
\end{eqnarray}
we can recast \eqref{eq4.8} in a compact form as
\begin{eqnarray}
\mathcal{J}(\underline{b})=\parallel Q\underline{b}-\underline{g}^{\epsilon} \parallel^{2}+\lambda \parallel \underline{b} \parallel^{2}, \label{eq4.10}
\end{eqnarray}
which is also known as the zeroth-order Tikhonov regularization. Higher-order regularization can also be employed by replacing the last term in \eqref{eq4.8} by the first-order derivative
\begin{eqnarray}
\parallel \underline{b}^{\prime}\parallel^{2}=\sum_{k=2}^{K}(b_{k}-b_{k-1})^{2},     \label{eq4.11}
\end{eqnarray}
or by the second-order derivative
\begin{eqnarray}
\parallel \underline{b}^{\prime \prime}\parallel^{2}=\sum_{k=3}^{K}(b_{k}-2b_{k-1}+b_{k-2})^{2},    \label{eq4.12}
\end{eqnarray}
etc. In general, $N\geq K$ and the minimization of \eqref{eq4.10} then yields the zeroth-order Tikhonov regularization solution
\begin{eqnarray}
\underline{b}_{\lambda}=(Q^{tr}Q+\lambda I)^{-1} Q^{tr}\underline{g}^{\epsilon}.  \label{eq4.13}
\end{eqnarray}

Once $\underline{b}$ has been found, the spacewise dependent force function is obtained using \eqref{eq4.2}. Also, the displacement solution $u(x,t)$ is obtained using \eqref{eq119} and \eqref{eq4.1}.

\section{Numerical Results and Discussion}
In this section, we illustrate and discuss the numerical results obtained using the combined BEM+Tikhonov regularization described in Sections 3 and 4.

For simplicity, we take $c=L=T=\mu=1$ and in the BEM we use constant time and space interpolation functions. We consider an analytical solution given by
\begin{eqnarray}
u(x,t)&=&\sin(\pi x)+t+\frac{t^{2}}{2},\ \ \ (x,t)\in[0,1]\times[0,\infty),  \label{eq5.1}\\
f(x)&=&1+\pi^{2}\sin(\pi x) \ \ \  x\in[0,1]. \label{eq5.2}  
\end{eqnarray}
This generates the input data \eqref{eq2}-\eqref{eq5} and \eqref{eq117} given by
\begin{eqnarray}
u(x,0)=u_{0}(x)=\sin(\pi x), \quad
u_{t}(x,0)=v_{0}(x)=1, \quad x\in[0,1],   \label{eq5.3}
\end{eqnarray}
\begin{eqnarray}
u(0,t)=p_{0}(t)=t+\frac{t^{2}}{2},\quad u(1,t)=p_{L}(t)=t+\frac{t^{2}}{2}, \quad t\in[0,\infty), \label{eq5.4}
\end{eqnarray}
\begin{eqnarray}
u_{x}(0,t)=q_{0}(t)=\pi , \quad t\in[0,1]. \label{eq5.5}
\end{eqnarray}

First we investigate  the performance of the BEM described in Section 3 to solve the direct well-posed problem \eqref{eq1110}-\eqref{eq14v} for the function $v(x,t)$. Remark that condition \eqref{eq3.3} is satisfied by the initial displacement $u_{0}(x)$ in \eqref{eq5.3} hence,  there is no need to employ the modified function \eqref{eq3.4}. Also note that the direct problem for $v(x,t)$ satisfying equation \eqref{eq1110} (with $c=1$) subject to the initial conditions \eqref{eq5.3} and the Dirichlet boundary conditions \eqref{eq5.4} does not have a closed form analytical solution available. 

Figure 2 shows the numerical results for $v_{x}(0,t)$, as a function of $t$, obtained using the BEM with various $M=N\in \lbrace 20,40,80 \rbrace$. From this figure it can be seen that a convergent numerical solution, independent of the mesh, is rapidly obtained. The numerical solution for $v_{x}(0,t)$ obtained at the points $(t_{n})_{n=\overline{1,N}}$ is then input into equation \eqref{eq4.4} to determine the values for $g(t_{n})$ and its noisy counterpart $g^{\epsilon}(t_{n})$ given by \eqref{eq4.7} for $n=\overline{1,N}$.
\begin{figure}[H]
\centering
\includegraphics[width=17cm]{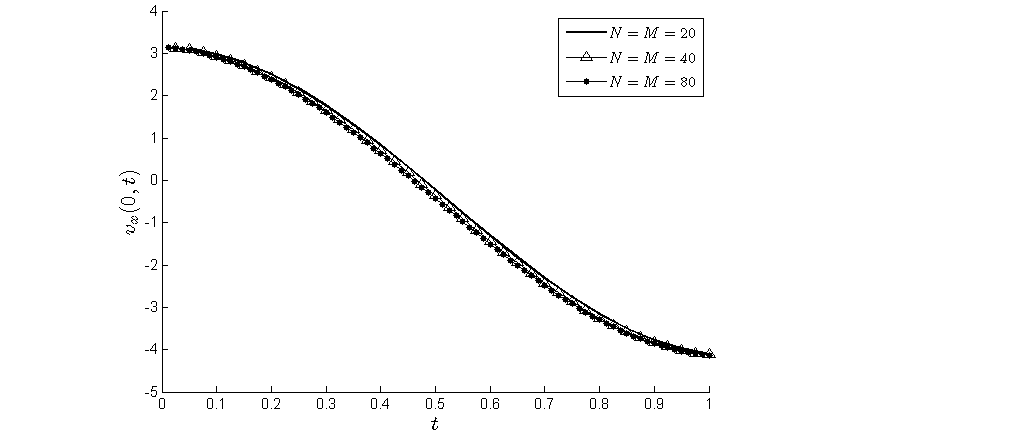}
\caption{The numerical results for $v_{x}(0,t)$ obtained using the BEM with $M=N\in \lbrace20,40,80 \rbrace$.}
\end{figure}
We turn now our attention to the pair solution \eqref{eq4.1} and \eqref{eq4.2} of the inverse problem \eqref{eq15w}-\eqref{eq19w}. Since this problem is ill-posed we expect that the matrix $Q$ in \eqref{eq4.9} having the entries 
\begin{eqnarray}
Q_{nk}=\frac{\sqrt{2}(1-\cos(\frac{k\pi n}{N}))}{k\pi}, \ \ \ n=\overline{1,N}, \ \ \ k=\overline{1,K}, \label{eq5.6} 
\end{eqnarray}
will be ill-conditioned. The condition number defined as the ratio between the largest to the smallest singular values of the matrix $Q$ is calculated in MATLAB using the command cond($Q$). Table 1 shows the condition number of the matrix $Q$ for various $N\in \lbrace 20,40,80 \rbrace$ and $K\in \lbrace 5,10,20 \rbrace$. We remark that the condition number is not affected by the increase in the number of measurements $N$, but it increases rapidly as the number $K$ of basis functions increases. The ill-conditioning nature of the matrix $Q$ can also be revealed by plotting its normalised singular values $sv(k)/sv(1)$ for $k=\overline{1,K}$, in Figure 3 for a fixed $N=80$ and $K=20$. These singular values have been calculated in MATLAB using the command svd($Q$).
\begin{figure}[H]
\centering
\includegraphics[width=17cm]{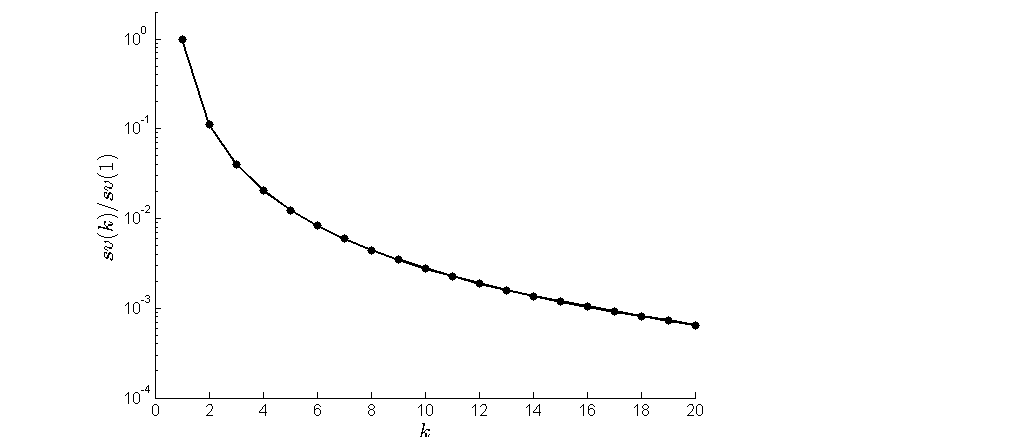}
\caption{Normalised singular values $sv(k)/sv(1)$ for $k=\overline{1,K}$, for $N=80$ and $K=20$.}
\end{figure}
\begin{table}
\caption{Condition number of the matrix $Q$ given by equation \eqref{eq5.6}.}
\centering
\begin{tabular}{|c|c|c|c|c|c}
\hline
$K$  &$N=20$   &$N=40$     &$N=80$ \\  
\hline
$5$	&$82.62$  &$82.25$  	&$82.28$	 \\
\hline
$10$	&$371.6$  &$367.0$  	&$365.7$	 \\
\hline
$20$	&$1.42E+3$  &$1.55E+3$  	&$1.54E+3$	 \\
\hline
\end{tabular}
\end{table}
Let us fix $N=80$ and now proceed to solving the inverse problem \eqref{eq15w}-\eqref{eq19w} which based on the method of Section 4 has been reduced to solving the linear, but ill-conditioned system of equations 
\begin{eqnarray}
Q\underline b=\underline g^{\epsilon}. \label{eq5.7}
\end{eqnarray} 
Using the Tikhonov regularization method one obtains a stable solution given explicitly by equation \eqref{eq4.13} provided that the regularization parameter $\lambda$ is suitably chosen. 
\subsection{Exact Data}
We first consider the case of exact data, i.e. $p=0$ and hence $\epsilon=0$ in \eqref{eq4.5} and \eqref{eq4.7}. Then $\underline{g}^{\epsilon}=\underline{g}$ and the system of equations \eqref{eq5.7} becomes  
 \begin{eqnarray}
Q\underline b=\underline g. \label{eq5.8}
\end{eqnarray} 
We remark that although we have no random noise added to the data $q_{0}$, we still have some numerical noise  in the data $g$ in \eqref{eq4.4}. This is given by the small discrepancy between the unavailable exact solution $v_{x}(0,t)$ of the direct problem and its numerical BEM solution obtained with $M=N=80$ plotted in Figure 2. However, the rapid convergent behaviour shown is Figure 2 indicates that this numerical noise is small (at least in comparison with the large amount of random noise $\epsilon$ that we will be including in the data $q_{0}$ in Section 5.2).

Figure 4 shows the retrieved coefficient vector $\underline b=(b_{k})_{k=\overline{1,K}}$ for $K=20$ obtained using no regularization, i.e. $\lambda =0$, in which case \eqref{eq4.13} produces the least-squares solution 
\begin{eqnarray}
\underline b=(Q^{tr}Q)^{-1}Q^{tr}\underline g \label{eq5.9}
\end{eqnarray} 
of the system of equations \eqref{eq5.7}.

Note that the analytical values for the sine Fourier series coefficients are given by
\begin{eqnarray*}
b_{k}=\sqrt{2}\int^{1}_{0}f(x)\sin(k\pi x)dx=\sqrt{2}\int^{1}_{0}(1+\pi^{2}\sin(\pi x))\sin(k\pi x)dx
\end{eqnarray*}
which gives
\begin{eqnarray}
b_{k} =
\begin{cases}
\frac{2\sqrt{2}}{\pi}+\frac{\pi^{2}}{\sqrt{2}}\simeq 7.8791 \ \ \ \ \ \ \ \ \ \quad \text{if} \ \ \ \ \ \ \ \ \ \ \quad k=1,
\\
0 \ \ \ \ \ \ \ \ \ \ \ \ \ \ \ \ \ \ \ \ \ \ \ \ \ \ \ \ \ \ \ \ \quad \text{if} \ \ \ \ \ \ \ \ \ \ \quad k=\text{even},
\\
\frac{2\sqrt{2}}{k\pi} \ \ \ \ \ \ \ \ \ \ \ \ \ \ \ \ \ \ \ \ \ \ \ \ \ \ \ \ \  \quad \text{if} \ \ \ \ \ \ \ \ \ \ \quad k=\text{odd}\geq 3.    \label{eq5.10}
\end{cases}
\end{eqnarray}

By inspecting Figure 4 it appears that the leading term $b_{1}$ is the most significant in the series expansions \eqref{eq4.1} and \eqref{eq4.2}. These expansions give the solutions $f(x)$ and $u(x,t)$ (via \eqref{eq119}) which are plotted in Figures 5 and 6, respectively. From these figures it can be seen that accurate numerical solutions are obtained.
\begin{figure}[H]
\centering
\includegraphics[width=17cm]{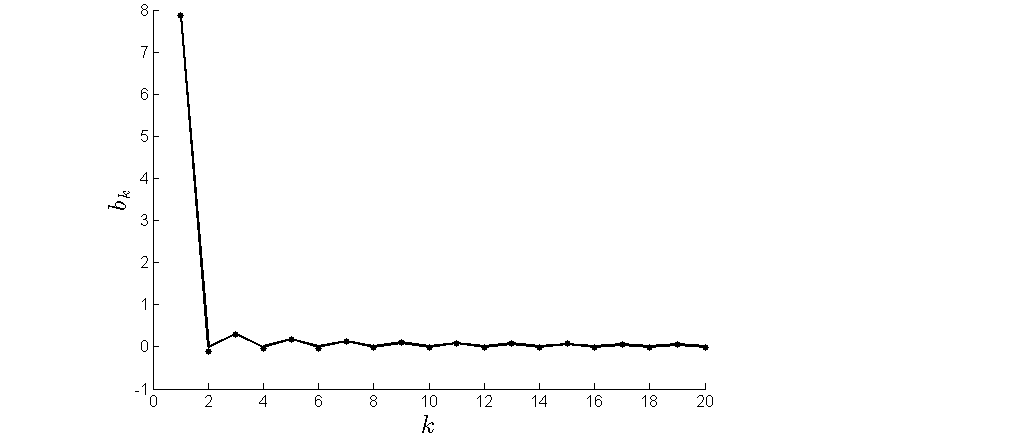}
\caption{The numerical solution ({\tt ...}) for $(b_{k})_{k=\overline{1,K}}$ for $K=20, \ N=80$, obtained with no regularization, i.e. $\lambda=0$, for exact data, in comparison with the exact solution \eqref{eq5.10} (-----).}
\end{figure}
\begin{figure}[H]
\centering
\includegraphics[width=17cm]{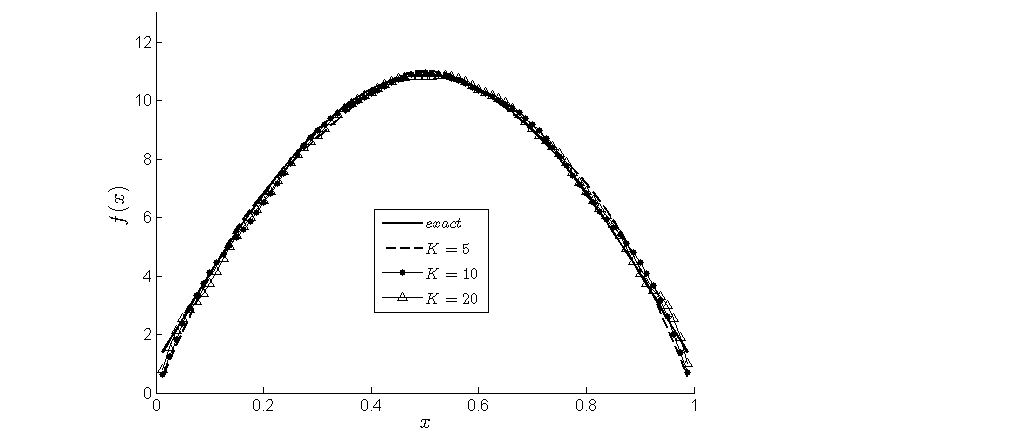}
\caption{The exact solution \eqref{eq5.2} for $f(x)$ in comparison with the numerical solution \eqref{eq4.2} for various $K\in \lbrace 5,10,20 \rbrace$, no regularization, for exact data.}
\end{figure}
\ \ \ \ \ \ \ \ \ 
\begin{figure}[H]
\includegraphics[width=8.5cm]{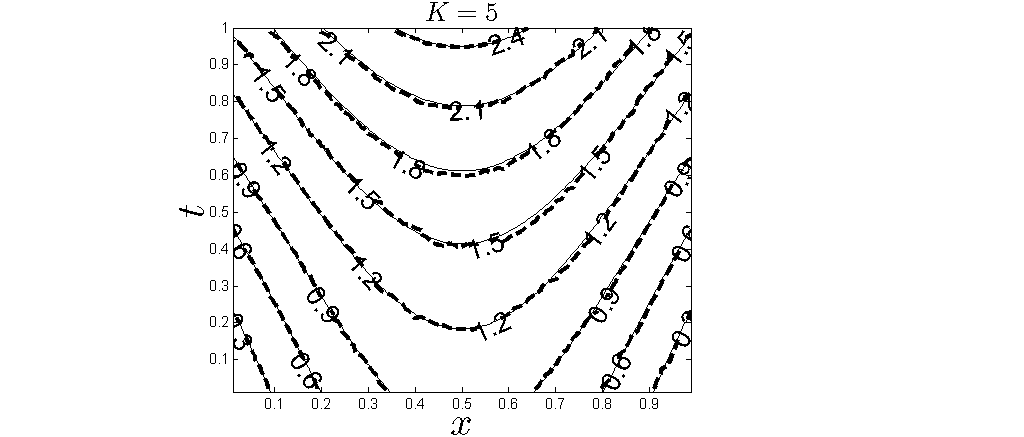}\includegraphics[width=8.5cm]{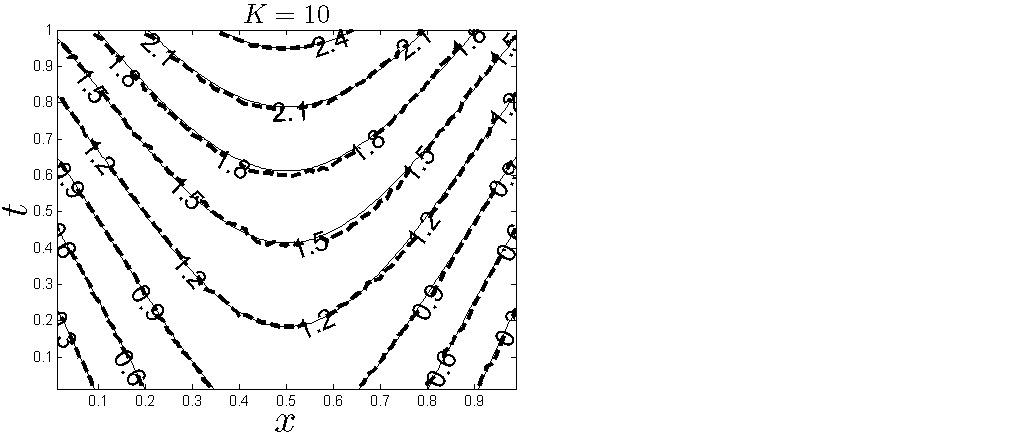}
\begin{center}
\includegraphics[width=8.5cm]{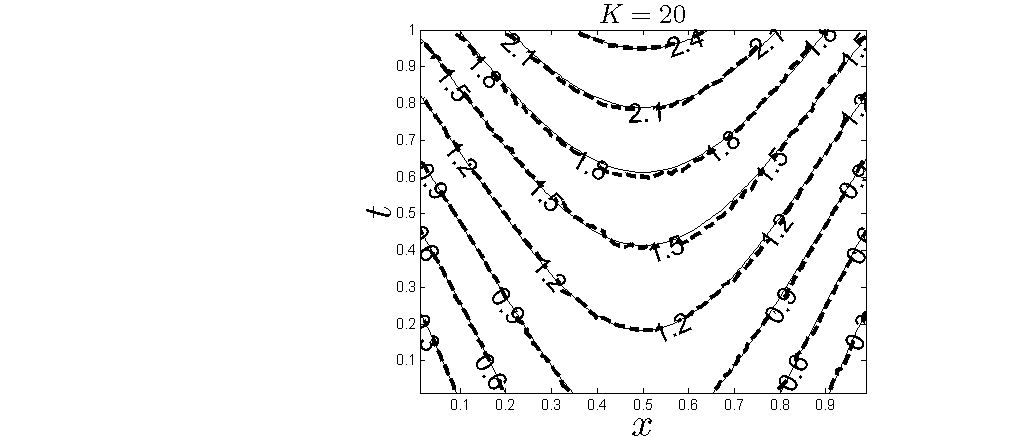}
\end{center}
\caption{The numerical solution $(---)$ for $u(x,t)$ obtained with various $K\in \lbrace 5,10,20 \rbrace$, no regularization, for exact data, in comparison with the exact solution \eqref{eq5.1} (-----).}
\end{figure}
\subsection{Noisy Data}
In order to investigate the stability of the numerical solution we include some ($p\%=1\%$) noise into the input data \eqref{eq117}, as given by equation \eqref{eq4.5}. The numerical solutions for $f(x)$ and $u(x,t)$ obtained for various values of $K\in \lbrace 5,10,20 \rbrace$ and no regularization are plotted in Figures 7 and 8, respectively. First, by inspecting Figures 6 and 8 it can be observed that there is little difference between the results for $u(x,t)$ obtained with and without noise and that there is very good agreement with the exact solution \eqref{eq5.1}. It also means that the numerical solution for the displacement $u(x,t)$ is stable with respect to noise added in the input data \eqref{eq117}. In contrast, in Figure 7 the unregularized numerical solution for $f(x)$ manifests instabilities as $K$ increases. For $K$ (small) the numerically retrieved solution is quite stable showing that taking a small number of basis functions in the series expansion \eqref{eq4.2} has a regularization effect. However, as $K$ increase to $10$ or $20$ it can be clearly seen that oscillations start to appear. 
\begin{figure}[H]
\centering
\includegraphics[width=17cm]{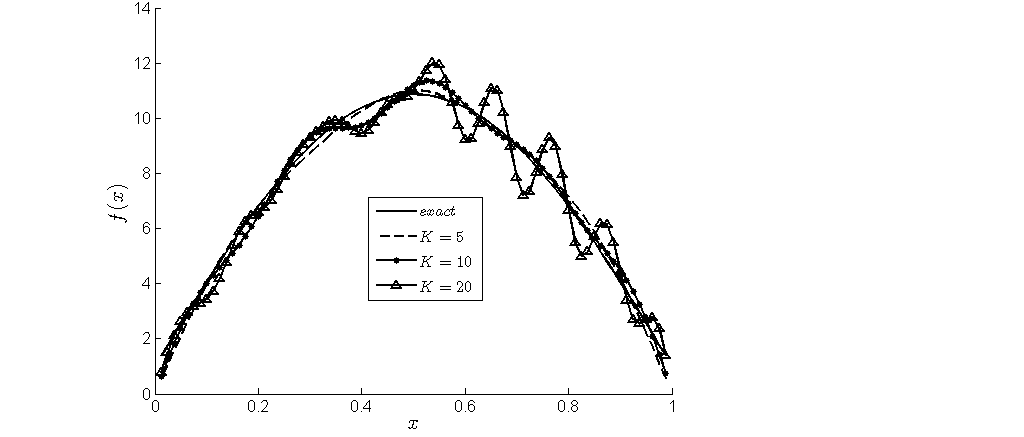}
\caption{The exact solution \eqref{eq5.2} for $f(x)$ in comparison with the numerical solution \eqref{eq4.2} for various $K\in \lbrace 5,10,20 \rbrace$, no regularization, for $p\%=1\%$ noisy data.}
\end{figure}
Eventually, these oscillations will become highly unbounded, as $K$ increases even further. In order to deal with this instability we employ the Tikhonov regularization which yields the solution \eqref{eq4.13}. We fix $K=20$ and we wish to alleviate the instability of the numerical solution for $f(x)$ shown by (---$\bigtriangleup$---) in Figure 7 obtained with no regularization, i.e. $\lambda=0$, for $p\%=1\%$ noisy data. Including regularization we obtain the solution \eqref{eq4.13} whose accuracy error, as a function of $\lambda$, is plotted in Figure 9. This error has been calculated as $||f_{numerical}-f_{exact}||=\sqrt{\sum_{n=1}^{N}(f_{numerical}(t_{n})-f_{exact}(t_{n}))^{2}}.$ From Figure 9 it can be seen that the minimum of the error occurs around $\lambda=10^{-1}$. Clearly, this argument cannot be used as a suitable choice for the regularization parameter $\lambda$ in the absence of an analytical (exact) solution \eqref{eq5.2} being available. However, one possible criterion for choosing $\lambda$ is given by the L-curve method, [6], which plots the residual norm $||Q\underline{b}_{\lambda}-\underline{g}^{\epsilon}||$ versus the solution norm $||\underline{b}_{\lambda}||$ for various values of $\lambda$. This is shown in Figure 10 for various values of $$\lambda \in \lbrace 10^{-3},5\times10^{-2},10^{-2},8\times10^{-2},6\times10^{-2},4\times10^{-2},2\times10^{-2},10^{-1},0.2,0.3,...,1\rbrace.$$ The portion to the right of the curve corresponds to large values of $\lambda$ which make the solution oversmooth, whilst the portion to the left of the curve corresponds to small values of $\lambda$ which make the solution undersmooth. The compromise is then achieved around the corner region of the L-curve where the aforementioned portions meet. Figure 10 shows that this corner region includes the values around $\lambda=10^{-1}$ which was previously found to be optimal from Figure 9.

Finally, Figure 11 shows the regularized numerical solution for $f(x)$ obtained with various values of the regularization parameter $\lambda \in \lbrace 10^{-2},10^{-1},10^{0}\rbrace$ for $p\%=1\%$ noisy data. From this figure it can be seen that the value of the regularization parameter $\lambda$ can also be chosen by trial and error. By plotting the numerical solution for various values of $\lambda$ we can infer when the instability starts to kick off. For example, in Figure 11, the value of $\lambda=10^{0}$ is too large and the solution is oversmooth, whilst the value of $\lambda=10^{-2}$ is too small and the solution is unstable. We could therefore inspect the value of $\lambda=10^{-1}$ and conclude that this is a reasonable choice of the regularization parameter which balances the smoothness with the instability of the solution.  
\begin{figure}[H]
\includegraphics[width=8.5cm]{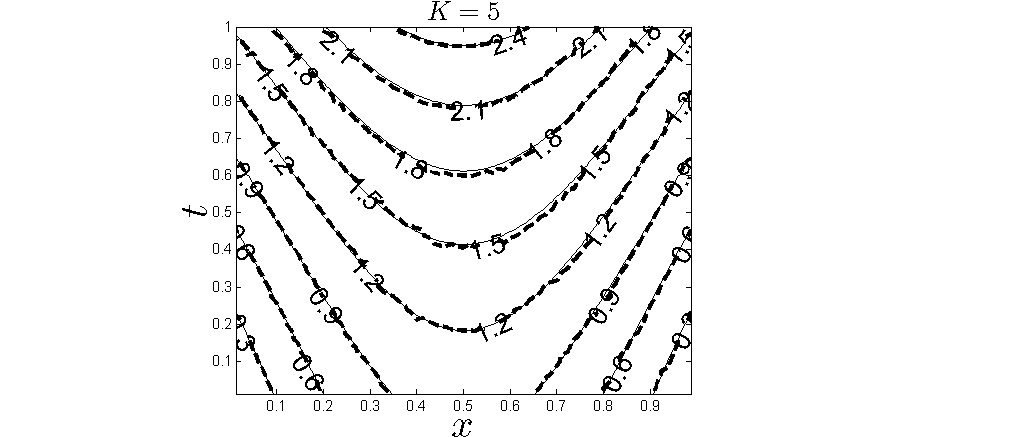}\includegraphics[width=8.5cm]{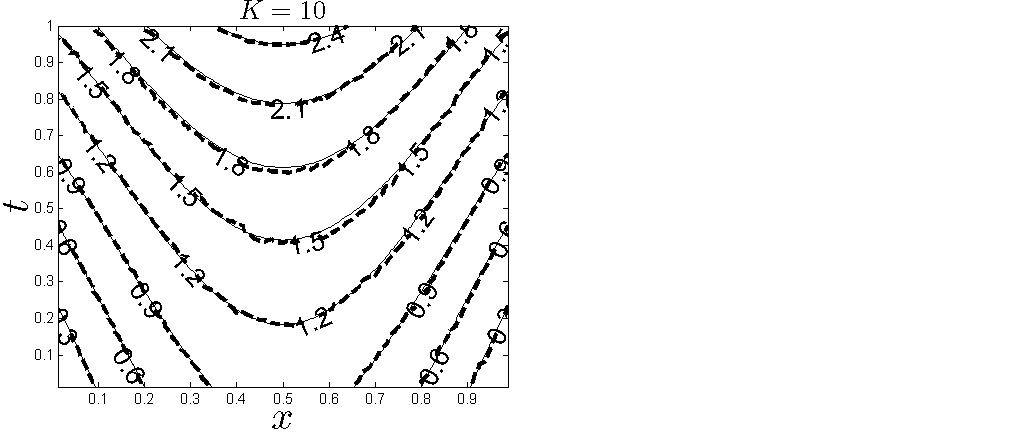}
\begin{center}
\includegraphics[width=8.5cm]{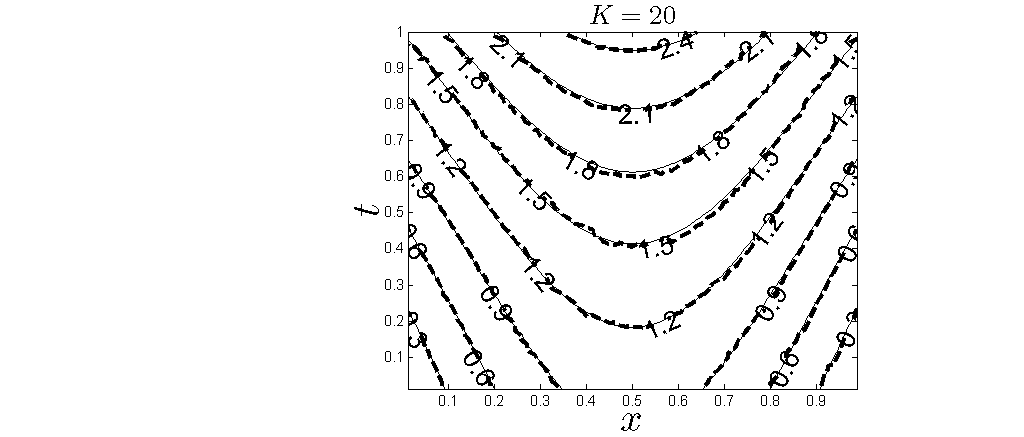}
\end{center}
\caption{The numerical solution $(---)$ for $u(x,t)$ obtained with various $K\in \lbrace 5,10,20 \rbrace$, no regularization, for $p\%=1\%$ noisy data, in comparison with the exact solution \eqref{eq5.1} (-----).}
\end{figure}
\begin{figure}[H]
\centering
\includegraphics[width=17cm]{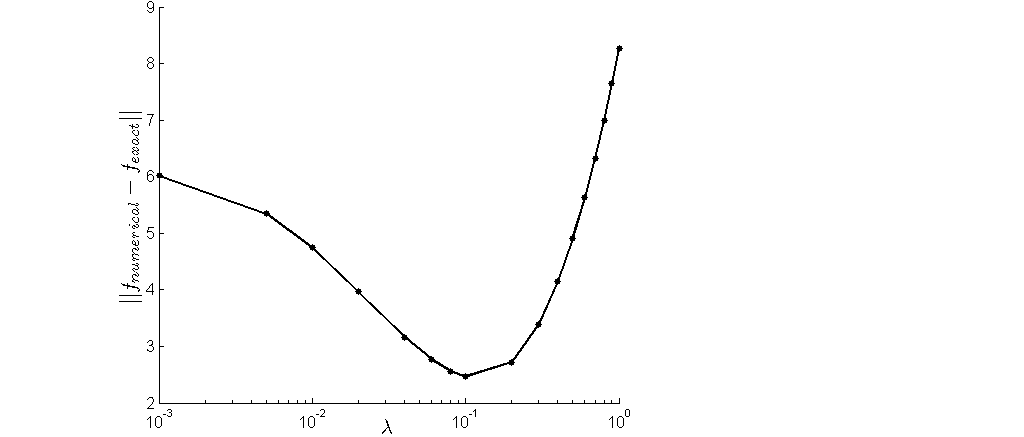}
\caption{The accuracy error $||f_{numerical}-f_{exact}||$, as a function of $\lambda$, for $K=20$ and $p\%=1\%$ noise. }
\end{figure}
\begin{figure}[H]
\centering
\includegraphics[width=17cm]{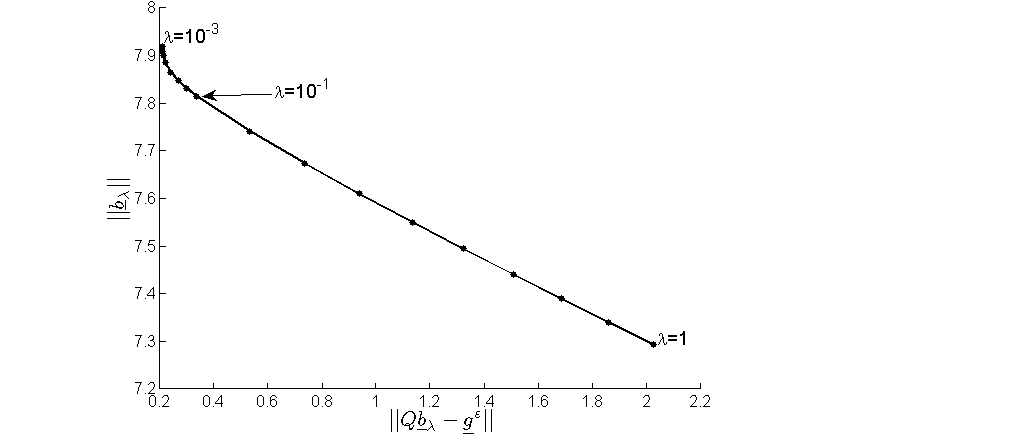}
\caption{The L-curve for the Tikhonov regularization  \eqref{eq4.10}, for $K=20$ and  $p\%=1\%$ noise.}
\end{figure} 
\begin{figure}[H]
\centering
\includegraphics[width=17cm]{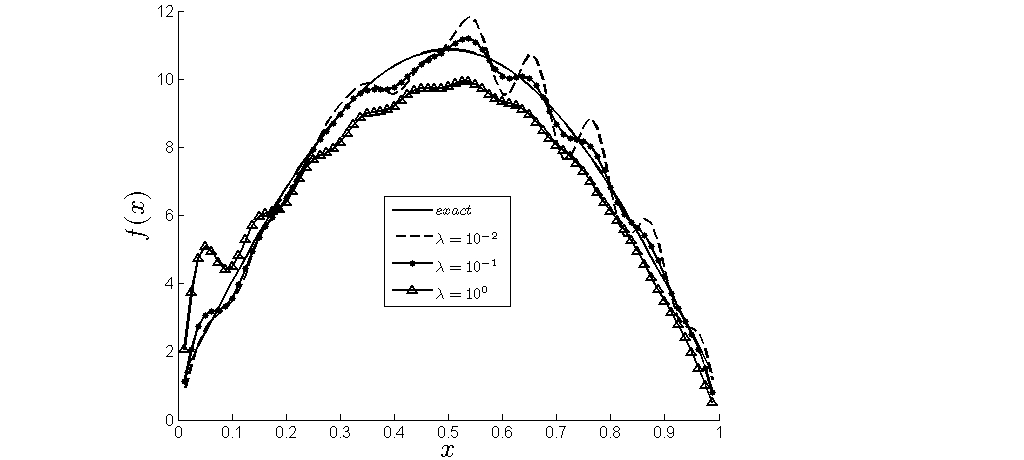}
\caption{The exact solution \eqref{eq5.2} for $f(x)$ in comparison with the numerical solution \eqref{eq4.2}, for $K=20$, $p\%=1\%$ noise, and regularization parameters $\lambda \in \lbrace 10^{-2},10^{-1},10^{0}\rbrace$.}
\end{figure}       
\section{Alternative Control}
For completeness, we describe the previously remarked, after equation \eqref{eq19w}, alternative control namely,
that we can replace equation \eqref{eq13v} by
\begin{eqnarray}
 v_{x}(0,t)=q_{0}(t), \quad t\in[0,\infty) \label{eq76}
\end{eqnarray}
and equation \eqref{eq19w} by
\begin{eqnarray}
 w(0,t)=p_{0}(t)-v(0,t), \quad t\in[0,T]. \label{eq77}
\end{eqnarray}
Then we can solve the well-posed direct problem \eqref{eq1110}-\eqref{eq12v}, \eqref{eq14v} and \eqref{eq76} to
obtain first $v(0,t)$. For the same test example, as in the previous section, Figure 12 shows the numerical
results for $v(0,t)$ obtained using the BEM with $M=N\in\lbrace20,40,80\rbrace$. From this figure it can be seen that a
convergent numerical solution, independent of the mesh, is rapidly achieved. The value of $v(0,t)$ is then
introduced into \eqref{eq77} to generate the Dirichlet data at $x=0$ for the inverse problem given by equation
\eqref{eq15w}-\eqref{eq18w} and \eqref{eq77}. We solve this inverse problem, as described in Section 4, with the 
obvious modifications to obtain the separation of variables solution
 \begin{eqnarray}
  w_{K}(x,t;\underline{b})=\frac{\sqrt{2}}{c^{2}}\sum_{k=1}^{K}\frac{b_{k}}{\lambda_{k}^{2}}(1-\cos(c\lambda_{k}t))
  \cos(\lambda_{k}x), \quad (x,t)\in[0,L]\times[0,\infty), \label{eq78}
 \end{eqnarray}
\begin{eqnarray}
 f_{K}(x)=\sqrt{2}\sum_{k=1}^{K}b_{k}\cos(\lambda_{k}x), \quad x\in(0,L), \label{eq79}
\end{eqnarray}
where $\lambda_{k}=(k-\frac{1}{2})\pi/L$ for $k=\overline{1,K}$. 
The coefficient $\underline{b}=(b_{k})_{k=\overline{1,K}}$ is determined by imposing the additional condition
\eqref{eq77},
\begin{eqnarray}
 p_{0}(t)-v(0,t)=:h(t)=\frac{\sqrt{2}}{c^{2}}\sum_{k=1}^{K}\frac{b_{k}}{\lambda_{k}^{2}}(1-\cos(c\lambda_{k}t)),
\quad t\in[0,T], \label{eq80}
\end{eqnarray}
in the Tikhonov regularized sense \eqref{eq4.8}, namely, as minimizing the functional 
\begin{eqnarray}
\mathcal{J}(\underline{b}):=\sum_{n=1}^{N}{\bigg[}\frac{\sqrt{2}}{c^{2}}\sum_{k=1}^{K}\frac{b_{k}}{\lambda_{k}^{2}}(1-\cos(c\lambda_{k} t_{n}))-h^{\epsilon}(t_{n}){\bigg]}^{2}+\lambda \sum_{k=1}^{K}b_{k}^{2}. \label{eq81}
\end{eqnarray}
Denoting
\begin{eqnarray}
h^{\epsilon}=(h^{\epsilon}(t_{n}))_{n=\overline{1,N}}, \quad
Q_{nk}=\frac{\sqrt{2}(1-\cos(c\lambda_{k} t_{n}))}{c^{2}\lambda_{k}^{2}}, \quad n=\overline{1,N}, \ k=\overline{1,K},
\label{eq82}
\end{eqnarray}
we can recast \eqref{eq81} in the compact form \eqref{eq4.9}.

The condition numbers of the matrix $Q$, defined in equation \eqref{eq82}, are given in Table 2 for various $N\in\lbrace20,40,80\rbrace$ and $K\in\lbrace5,10,20\rbrace$. From this table it can be seen that ill-conditioning increases significantly, as $K$ increases.
\begin{figure}[H]
\centering
\includegraphics[width=17cm]{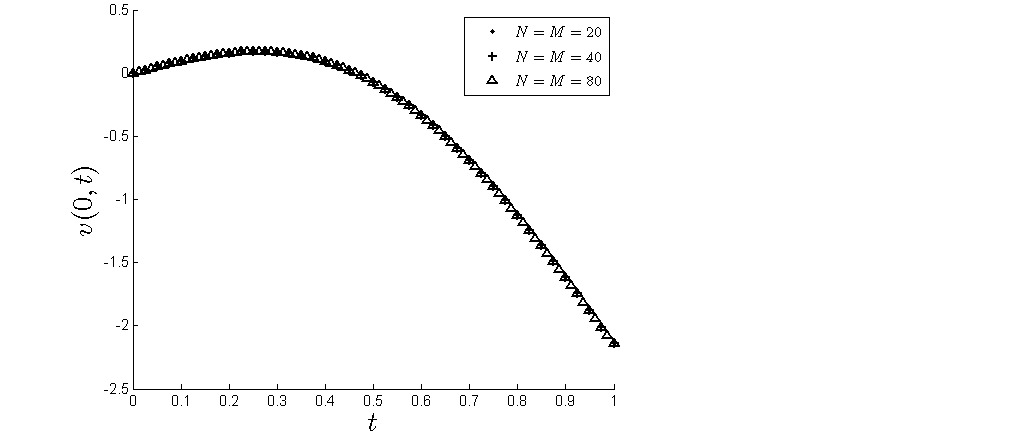}
\caption{The numerical results for $v(0,t)$ obtained using the BEM with $M=N\in \lbrace20,40,80 \rbrace$.}
\end{figure} 
\begin{table}
\caption{Condition number of the matrix $Q$ given by equation \eqref{eq82}.}
\centering
\begin{tabular}{|c|c|c|c|c|c}
\hline
$K$  &$N=20$   &$N=40$     &$N=80$ \\  
\hline
$5$	&$ 3.55E+3$  &$3.62E+3$  	&$3.68E+3$	 \\
\hline
$10$	&$6.81E+4$  &$6.84E+4$  	&$6.96E+4$	 \\
\hline
$20$	&$1.21E+6$  &$1.17E+6$  	&$1.18E+6$	 \\
\hline
\end{tabular}
\end{table}
\subsection{Exact Data}
In the case of exact data, as in Section 5.1, Figure 13 shows the retrieved coefficients $(b_{k})_{k=\overline{1,K}}$ for $K=20$ obtained with no regularization, i.e. $\lambda=0$, in comparison with the exact cosine Fourier series coefficients given by
\begin{eqnarray*}
b_{k}=\sqrt{2}\int_{0}^{1}f(x)\cos\left(\left(k-\frac{1}{2}\right)\pi x\right)dx
\end{eqnarray*}
which, for $f(x)$ given by \eqref{eq5.2}, gives
\begin{eqnarray}
b_{k} =
\begin{cases}
\frac{2\sqrt{2}(2\pi^{2}+3)}{3\pi}\simeq 6.8242 \ \ \ \ \ \ \ \ \ \ \ \ \quad \quad \quad \quad \quad \quad \text{if} \ \ \ \ \ \ \ \ \ \quad k=1,
\\
\\
-\frac{2\sqrt{2}(2\pi^{2}(2k-1)+(-1)^{k}(4k^{2}-4k-3))}{\pi(8k^{3}-12k^{2}-2k+3)} \ \ \ \ \ \ \ \ \ \ \ \ \ \quad \text{if} \ \ \ \ \ \ \ \ \ \quad k>1. \label{eq83}
\end{cases}
\end{eqnarray}
Good agreement between the exact and numerical values can be observed. With these value of $\underline{b}=(b_{k})_{k=\overline{1,K}}$, the solution \eqref{eq79} for the force function yields the numerical results illustrated in Figure 14. From this figure it can be seen that accurate numerical results are obtained.
\begin{figure}[H]
\centering
\includegraphics[width=17cm]{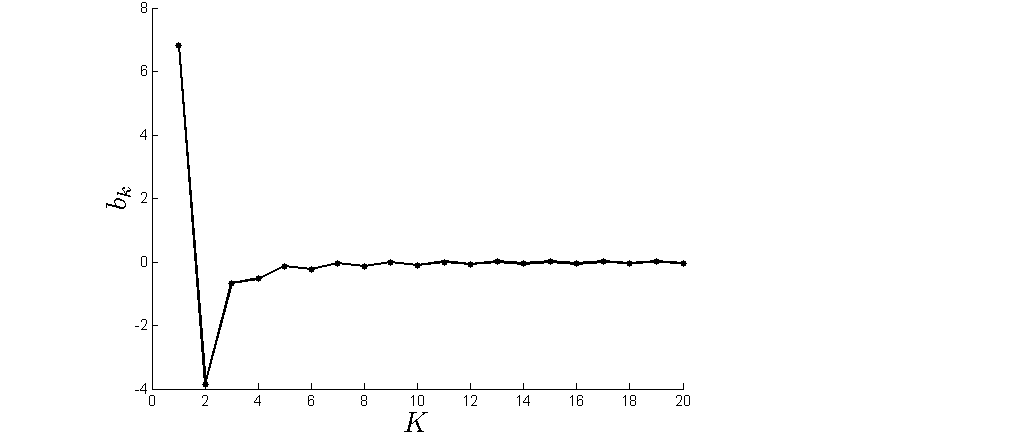}
\caption{The numerical solution ({\tt ...}) for $(b_{k})_{k=\overline{1,K}}$ for $K=20, \ N=80$, obtained with no regularization, i.e. $\lambda=0$, for exact data, in comparison with the exact solution \eqref{eq83} (-----).}
\end{figure} 
\begin{figure}[H]
\centering
\includegraphics[width=17cm]{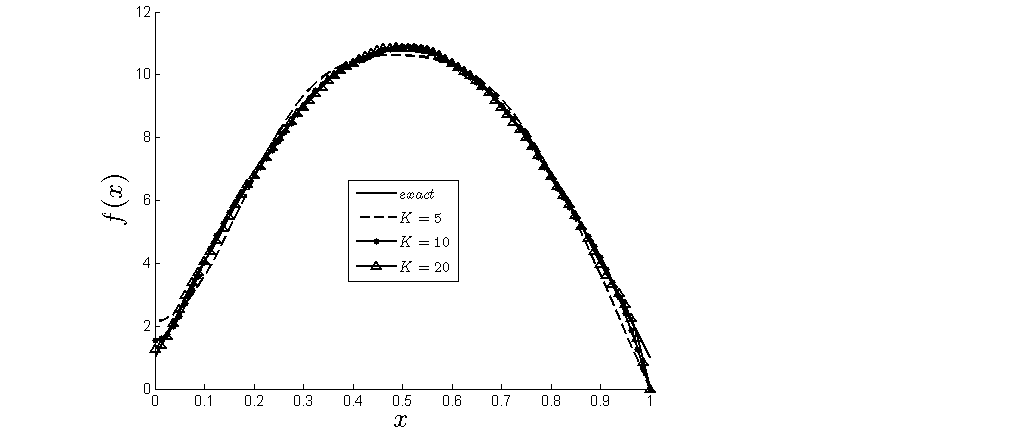}
\caption{The exact solution \eqref{eq5.2} for $f(x)$ in comparison with the numerical solution \eqref{eq79} for various $K\in \lbrace 5,10,20 \rbrace$, no regularization, for exact data.}
\end{figure} 
\subsection{Noisy Data}
In the case of noisy data, as in Section 5.2, Figure 15 shows the regularized numerical solution for $f(x)$ obtained with various $\lambda \in\lbrace 10^{-4},10^{-3},10^{-2}\rbrace$ for $p\%=1\%$ noisy data added to $p_{0}(t)$ as
\begin{eqnarray}
p_{0}^{\epsilon}(t_{n})=p_{0}(t_{n})+\epsilon, \quad n=\overline{1,N}, \label{eq84}
\end{eqnarray} 
where $\epsilon$ are $N$ random noisy variables generated from a Gaussian normal distribution with mean zero and standard deviation $\sigma$ given by
\begin{eqnarray}
\sigma=p\%\times max_{t\in[0,T]}\left|p_{0}(t)\right|. \label{eq85} 
\end{eqnarray}
From this figure it can be seen that the numerical results obtained with $\lambda$ between $10^{-3}$ and $10^{-2}$ are reasonably stable and accurate.
\begin{figure}[H]
\centering
\includegraphics[width=17cm]{{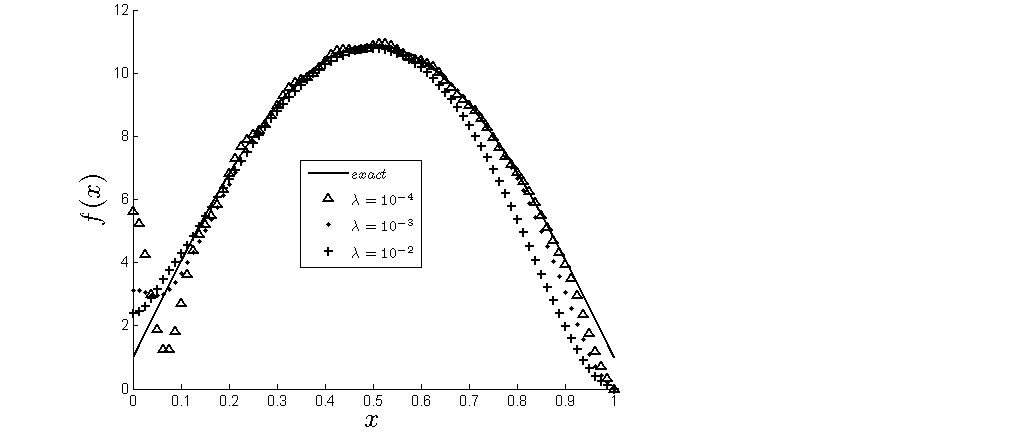}}
\caption{The exact solution \eqref{eq5.2} for $f(x)$ in comparison with the numerical solution \eqref{eq79} for $K=20$, $p\%=1\%$ noise, and regularization parameters $\lambda \in \lbrace 10^{-4},10^{-3},10^{-2}\rbrace$.}
\end{figure}       
\section{Conclusions}
An inverse force problem for the one-dimensional wave equation has been investigated. The unknown forcing term was assumed to depend on the space variable only and the additional measurement which ensures a unique retrieval was the flux at one end of the string. This inverse problem is uniquely solvable, but is still ill-posed since small errors in the input flux cause large errors in the output force. The problem is split into a direct well-posed problem for the linear wave equation, which is solved numerically using the BEM, and an inverse ill-posed problem whose unstable solution is expressed as a separation of variables truncated series. In order to stabilise the solution, the Tikhonov regularization method has been employed. The choice of the regularization parameter was based on the L-curve criterion. Numerical results for a typical benchmark smooth test example show that an accurate and stable solution has been obtained. Future work will consist in investigating the inverse force problem in higher dimensions.\\
\\  
\textbf{\large Acknowledgment}\\
S.O. Hussein would like to thank the Human Capacity Development Programme (HCDP) in Kurdistan for their financial support in this research.

\newpage


\newpage


\end{document}